\documentclass[12pt]{article}

\newtheorem{prop}{Proposition}
\newtheorem{thm}{Theorem}

\newtheorem{cor}{Corollary}
\newcommand{\Nabla }{\nabla}

\newcommand{\doublearrow}{\begin{picture}(15,12)(0,3)
\put(0,0){$\to$}
\put(0,8){$\to$}
\end{picture}}
\title{Algebra of Principal Fibre Bundles, and Connections.}
\author{Anders Kock}

\begin{document}
\maketitle

In this paper, I intend to put together some of the efforts by several 
people  of making 
aspects of fibre bundle theory into algebra.
The initiator of these efforts was Charles Ehresmann, who put the notion of 
{\em groupoid}, and {\em groupoid action} in the focus for fibre 
bundle theory in general, and for connection theory in particular.

In so far as connection theory is concerned, this paper is a 
sequel to \cite{CCBI}, and we presuppose some of the notions 
presented there: those of Sections 1, 3, 7, 8, and 11, so they 
will be recalled only sketchily. (The paper may also partly be seen 
as a rewiting of \cite{CN}.)


\section {Principal Fibre Bundles} Let us consider a  groupoid object $\Phi \doublearrow C$ in a  left exact 
category $\underline{E}$. Let us also consider a 
subobject $A\subseteq C$ and a global section $*:1\to C$.
 We shall talk about $\underline{E}$ as if 
it were the category of sets, so we may say ``subset'' instead of ``subobject''. 
We assume that the domain- and codomain formation maps are effective 
descent maps in $\underline{E}$, and that the groupoid is {\em 
transitive}, meaning that ``the anchor'' map $<d_0 ,d_1 > :\Phi \to 
C\times C$ is also an effective descent map.

Then the set $P$ of those arrows of $\Phi$, whose codomain is in $A$ and whose domain 
is $*$, carries the structure of a principal fibre bundle over $A$, 
with group $G=\Phi (*,*)$. 
 Any principal fibre bundle in 
$\underline{E}$ comes about this way from a groupoid (see the remarks below). 
The algebraic structure of $P$ 
comes about from that of $\Phi$, and may be made explicit as follows. 
First, codomain formation $d_1 :\Phi \to C$ resticts to a map $\pi :P\to 
A$, which is the structural map of the bundle. The group $G=\Phi (*,*)$ acts from the 
right on $P$, by precomposition (we compose from right to left). 
Clearly, this action is free, 
and transitive on the fibres $\pi ^{-1}(a)$ ($a\in A$). 

Any element $g$ of $G=\Phi (*,*)$ may, for any $a\in A$, be written in the form 
$x^{-1}z$ for a pair of elements in $\pi ^{-1}(a)$. This 
representation of elements in $G$ by ``fractions'' $x^{-1}z$ prompts us to 
use the (Ehresmann) notation $P^{-1}P$  for $G$. Then clearly $x\cdot 
x^{-1}z =z$ (where $\cdot$ denotes the $G$-action). Any choice of $x 
\in \pi ^{-1} (a)$ provides us with an explicit bijection $ \pi ^{-1} 
(a) \to G$, given by $z \mapsto x^{-1}z$. 

Let us also consider the (transitive) subgroupoid $\Phi _A$ of $\Phi$ consisting 
of those arrows whose domain and codomain both belong to $A$. The 
groupoid $\Phi _A$ acts on the left on $P\to A$ by postcomposition; any arrow 
$a\to b$ in 
it may be presented in the form $y x^{-1}$, for some $x\in \pi 
^{-1}(a)$ and $y\in \pi ^{-1} (b)$. Then clearly $y x^{-1} \cdot 
x = y$ (where $\cdot $ denotes the action). The representation of 
arrows in $\Phi _A$ by ``fractions'' $y x^{-1}$ prompts us to use 
the (Ehresmann) notation $PP^{-1}$ for $\Phi _A$. 

\medskip

{\bf Remark.} The set $P$ itself carries a partially defined ternary 
operation, given by the composite $yx^{-1}z$ in $\Phi$ (defined 
subject to the book-keeping condition that $\pi (x)= \pi (z)$), and 
this operation satisfies a couple of equations and book-keeping 
conditions, making it into a ``pregroupoid'' on $A$, in the sense of 
\cite{CN}. Out of such 
pregroupoid, a transitive groupoid $\Phi$ on $A+1$ may be 
constructed, which in turn gives rise to $P$ by the procedure 
described above (provided $\pi :P \to A$ is an effective descent 
map); this is in essence demonstrated in \cite{CN}. Principal 
fibre bundles (in the classical sense) $P \to A$,
 in the category of smooth manifolds, say, may, in a rather evident 
way, be provided with pregroupoid structure. So our
``groupoid theoretic'' way of describing the notion of principal 
fibre bundle subsumes the classical notion, and it is essentially 
Ehresmann's conception. --- A (non-transitive) generalization where $*:1 \to C$ is 
replaced by a  subset $B\subseteq C$, is considered in 
\cite{GFB}; this generalization is relevant for foliation 
theory, cf. loc.cit.\ and \cite{MAF}.  

\medskip

We shall henceforth be interested  in the case where the 
``base'' $A$ of the bundle $P\to A$ is to be thought of as a 
manifold, so we denote it by $M$ rather than by $A$.

\medskip

{\bf Remark} on fibre bundles in general. A principal fibre 
bundle $P\to M$  
with group $G$, may by the above be identified 
with a groupoid $\Phi$ with set of objects $M+1$, (and 
with $G=\Phi 
(*,*)$, where $*$ is the isolated point of $M+1$).  Similarly, a
fibre bundle $\pi :E\to M$, with  associated principal bundle $P$ 
and with fibre a left $G$-set $ F$, becomes identified with  a discrete 
opfibration  over $\Phi$ (in the algebraic sense, i.e.\ an {\em action} by 
$\Phi$), with $F = \pi ^{-1} (*)$, and $E=\pi 
^{-1}(M)$. Such fibre bundle is determined up to isomorphism by $P$ 
and $F$ (with its left $G$-action). In the present general context, 
this is the upshot of \cite{FBGC}. We shall not 
explicitly be using this correspondence for general fibre bundles 
here. But let us remark that $P$ itself is a fibre bundle, with 
fibre $G$ (with $G$-action by  left multiplication). Likewise, if we 
let $G$ act on $G$ by conjugation, $g\cdot h := ghg^{-1}$, we get a 
group bundle, namely what \cite{Mackenzie} calls the {\em  
gauge group bundle}  of $PP^{-1}$. (It is also known under the 
name $Ad(P)$.) 
We shall utilize this latter bundle, but shall recall it without 
reference to this general fibre bundle theory. For a groupoid $\Psi 
\doublearrow M$, the gauge group bundle $\mbox{{\bf gauge}}(\Psi )$ is a bundle over $M$, 
which  for its fibre over $a\in 
M$ simply has the group $$(\mbox{{\bf gauge}}(\Psi ))_a = \Psi (a,a).$$ 
It carries a left action by $\Psi$, 
given by conjugation: if $ f: a \to b$ in $\Psi$ and $h \in  \Psi 
(a,a)$, then $fhf^{-1} \in \Psi (b,b)$. 
\medskip

The existence, for any principal fibre bundle $P$,    of an embedding
of $P$ into a groupoid $\Phi$,   
 implies a ``metatheorem'', namely 
that we may calculate freely with expressions, like $vu^{-1}$, 
as if we were dealing with actual compositions in a groupoid.
The 'action' dots, like in $yx^{-1}\cdot x$ are then 
superfluous, and the same applies to many parentheses; so they 
are mainly kept for readability.
 The message (which I also 
tried to get through in \cite{CCBI} and in several other places) 
is that a fair amount of calculations 
in geometry can be performed on this very basic ``multiplicative'' 
level.

\medskip

Since an arrow $f:a\to b$ in the groupoid $PP^{-1}$ may be represented 
as a ``fraction'' $yx^{-1}$ (with $y\in P_b$ and $x\in P_a$), it follows 
that an element $h$ over $a$ in the 
gauge group bundle $\mbox{gauge}(PP^{-1})$ may be represented by a 
fraction $yx^{-1}$ with $y$ and $x$ both $\in P_a$.
For the case where the group $G$ is {\em commutative}, it is well known, and 
easy to see, that we have an isormorhism of group bundles
\begin{equation} 
\mbox{{\bf gauge}}(PP^{-1}) \cong M\times G,
\label{comm}\end{equation}
given by sending $h=yx^{-1}\in PP^{-1}$ to $x^{-1}y \in P^{-1}P$. 
This cannot be done for non-commutative $G$: for any $g\in G$, the same $h$ may also 
be represented by the fraction 	$yg(xg)^{-1}$, but 
$(xg)^{-1}(yg)=g^{-1}(x^{-1}y)g$ which is not equal to $x^{-1}y$ in 
general.

\section{Connections versus connection forms}

Consider a principal bundle $\pi : P \to M$, with group $G$, as above.
 We shall  assume that $M$ and $P$ are equipped with  reflexive symmetric 
relations $\sim$, called the {\em neighbour} relation. The set of pairs $(x,y)\in M\times M$
 with $x\sim y$ is a subset 
$M_{(1)} \subseteq M\times M$, called the {\em first neighbourhood of 
the diagonal}, and similarly for $P_{(1)}\subseteq P\times P$. 
 We assume that $\pi :P\to M$ preserves the 
relation $\sim$, and also that it is an ``open submersion'' in the 
sense that if $a\sim b$ in $M$, and $\pi (x)=a$, then there exists a 
$y\sim x$ in $P$ with $\pi (y)=b$. In fact, we assume that for any 
``infinitesimal $k$-simplex'' $a_0 ,\ldots ,a_k$ in $M$ (meaning a 
$k+1$-tuple of mutual neighbours), and for any $x_0 \in P$ above $a_0$, 
there exists an infinitesimal $k$-simplex $x_0 ,\ldots ,x_k$ in $P$ 
(with the given first vertex $x_0$) which by $\pi $ maps to $a_0 
,\ldots ,a_k$.
 Finally. the action of any $g\in G$ 
on $P$ is assumed to preserve the relation $\sim$ on $P$.

This is motivated by Synthetic Differential Geometry (SDG), 
cf.\ \cite{SDG}, and more recently \cite{CCBI}, where the notion of 
connection (infinitesimal parallel transport) and differential form
 is elaborated in these 
terms. 
 
The groupoid viewpoint for connections is also in essence due to 
Ehresmann. In SDG, this connection notion becomes paraphrased (see 
\cite{CN}, \cite{ILG} or \cite{CCBI}, Section 8):
for a groupoid $\Phi \doublearrow M$, a connection in it is 
just a map $\nabla :M_{(1)} \to \Phi$ of reflexive symmetric graphs over $M$.

Let $\pi :P\to M$ be a principal fibre bundle. To any 
connection $\Nabla $ in the groupoid $PP^{-1}$, one may 
associate a 1-form $\omega $ on $P$ with values in the group $P^{-1}P$, 
as follows. For $u$ and $v$ neighbours in $P$, with $\pi 
(u)=a$, $\pi (v)= b$, put
\begin{equation}
\omega (u,v) := u^{-1}(\Nabla (a,b)\cdot v).
\label{cf}\end{equation}
Note that both $u$ and $\Nabla (a,b)\cdot v$ are in the $\pi 
$-fibre over $a$, so that the "fraction" $u^{-1}(\Nabla 
(a,b)\cdot v)$ makes sense as an element of $P^{-1}P$.

The defining equation is equivalent to
\begin{equation}u\cdot \underbrace{\omega(u,v)}_{\in P^{-1}P} = 
\underbrace{\nabla (\pi (u),\pi (v))}_{\in PP^{-1}}\cdot v .
\label{basic}\end{equation}
If  we agree that (for $u ,v$ in $P$ a pair of neighbours in $P$) 
 $\nabla 
(u,v)$ denotes $\nabla (\pi (u), \pi 
(v))$, this equation may be written more succinctly
\begin{equation}u\cdot\omega (u,v) = \nabla (u,v)\cdot v.
\label{c3}\end{equation}

It is possible to represent the relationship between $\nabla$ and the 
associated $\omega$ by means of a simple figure:

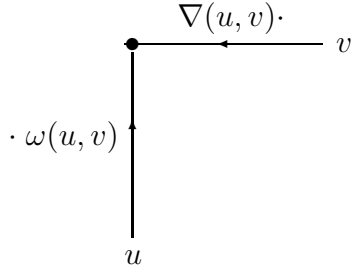
\begin{center}
\begin{picture}(100,100)
\put(0,0){$u$}
\put(80,80){$v$}
\put(0,80){$\bullet$}
\put(0,83){\line(1,0){75}}
\put(40,83){\vector(-1,0){5}}
\put(3,10){\line(0,1){70}}
\put(3,50){\vector(0,1){5}}
\put(-44,45){$\cdot \; \omega (u,v)$}
\put(20,90){$\nabla (u,v)\cdot $}
\end{picture}
\end{center}
The figure reflects something geometric, namely that $\omega (u,v)$ 
acts inside the fibre (vertically), whereas $\nabla$ defines a notion 
of horizontality.

We have the following two equations for $\omega$. First, let $x\sim 
y$ in $P$, and assume that $g$ has the property that also $xg\sim y$. 
Then
\begin{equation}
\omega (xg,y) = g^{-1}\omega (x,y).
\label{c1}\end{equation}
Also, for $x\sim y$ and {\em any} $g\in G$
\begin{equation}
\omega (xg,yg) = g^{-1}\omega (x,y) g.
\label{c2}\end{equation}
To prove (\ref{c1}), let us denote $\pi (x)= \pi (xg)$ by $a$ and 
$\pi (y)$ by $b$. 
Then we have, using the defining equation (\ref{basic}) for $\omega$ twice,
$$xg\;\omega (xg ,y) = \nabla(a,b)y = x\omega (x,y),$$
and now we may calculate as in a groupoid: first cancel the $x$ on 
the left, then multiply the equation  by $g^{-1}$ on the left.
To prove (\ref{c2}),
we have, with $a$ and $b$ as above, 
$$xg\; \omega(xg,yg) = \nabla (a,b) yg = x\omega(x,y)g,$$
by the defining equation (\ref{basic}) for $\omega (xg,yg)$, and
by (\ref{basic}) for $\omega (x,y)$, multiplied on the right by $g$,
respectively.
From this,  we get the result by first cancelling $x$ and then multiplying 
the equation by $g^{-1}$ on the left.

\medskip

The following Proposition is now the rendering, in our 
context, of the relationship between a connection $\nabla$ and 
its connection 1-form $\omega$:

\begin{prop}
The process $\nabla \mapsto \omega$ just described, establishes a 
bijective corresondence between 1-forms $\omega$ on $P$, with values in 
the group $P^{-1}P$ and
satisfying (\ref{c1}) and (\ref{c2}), and connections $\nabla$ in the groupoid 
$PP^{-1}$.
\label{one}\end{prop}
{\bf Proof.} Given a 1-form $\omega$ satisfying (\ref{c1}) and 
(\ref{c2}), we construct a connection $\nabla$ as follows. Let $a\sim 
b$ in $M$. To define the arrow $\nabla (a,b)$ in $PP^{-1}$, pick 
$u\sim v$ above $a\sim b$, and put
$$\nabla (a,b) =u(v\omega(v,u))^{-1}.$$
We first argue that this is independent of the choice of $v$, once $u$ 
is chosen. Replacing $v$ by $vg\sim u$, we are in the situation 
where (\ref{c1}) may be applied; we get
$$u(vg\; \omega (vg,u)^){-1}= u (vg\; g^{-1}\omega(v,u))^{-1}= 
u(v\omega(v,u))^{-1};$$
 the left hand side is $\nabla (a,b)$ defined using $u,vg$, the 
right hand side is using $u,v$.

To prove independence of choice of $u$: any other choice is of form 
$ug$ for some $g\in G$. For our new $v$, we now chose $vg$ (the 
result will not depend on the choice, by the argument just given).
Again we calculate. By (\ref{c2}), we have the first equality sign in
$$ug(vg\omega(vg,ug))^{-1}= ug(vg g^{-1}\omega (u,v)g)^{-1} =
ug(v\omega(u,v)g)^{-1} =u(v\omega (u,v))^{-1},$$
and the two expressions here are $\nabla (a,b)$ defined using, 
respectively, $ug,vg$ and $u,v$.

The calculation that the two processes are inverse of each other is  
trivial (using $\omega (u,v)=\omega (v,u)^{-1}$ and $\nabla 
(a,b)=\nabla (b,a)^{-1}$).

\section{Gauge forms versus horizontal equivariant forms}
We consider a principal fibre bundle $\pi :P\to M$ as in the 
previous section.  The {\em  horizontal} $k$-forms that we now consider, 
are $k$-forms on $P$ with values in the group $G =P^{-1}P$. 
{\em Horizontality} means for 
a $k$-form $\theta$ that
\begin{equation}
\theta (u_0 ,u_1 ,\ldots ,u_k )=\theta (u_0  ,u_1 \cdot g_1 ,\ldots 
,u_k \cdot g_k )
\label{bb} \end{equation}
for any infinitesimal $k$-simplex $ (u_0 ,u_1 ,\ldots ,u_k )$ 
in $P$, and any $g_1 ,\ldots g_k\in P^{-1}P$ with the property that
$(u_0  ,u_1 \cdot g_1 ,\ldots 
,u_k \cdot g_k )$ is still an infinitesimal 
simplex (which is a strong "smallness" requirement on the $g_i$'s). 

Note that the connection form $\omega$ for a connection $\nabla$ is 
{\em not} a horizontal 1-form, since $\omega (x,yg)=\omega(x,y) 
g$, not $=\omega (x,y)$.

We say that a $k$ form $\theta $, as above, is {\em equivariant} 
if for any infinitesimal $k$-simplex $(u_0 ,\ldots ,u_ k )$, 
and {\em any} $g\in P^{-1}P$, we have
\begin{equation}
\theta (u_0 \cdot g ,u_1 \cdot g ,\ldots ,u_k \cdot g ) = 
g^{-1} \theta (u_0 ,u_1 ,\ldots ,u_k ) g.
\label{ccc} \end{equation}
Note that connection forms are equivariant in this sense, by (\ref{c2}).

\begin{prop} Assume that the group $G=P^{-1}P$ is {\em commutative}. 
Then 
any horizontal equivariant $k$-form $\theta$ on $P$ can be written 
$\pi ^* (\Theta )$ for a unique $G$-valued $k$-form $\Theta$ on the 
base space $M$.
\label{two}\end{prop}

{\bf Proof.} It is evident that any form $\pi ^* (\Theta )$ is 
horizontal and equivariant (which here is better called {\em 
invariant}, since the equivariance condition now reads
$ \theta (u_0 \cdot g ,u_1 \cdot g ,\ldots ,u_k \cdot g ) = 
\theta (u_0 ,u_1 ,\ldots ,u_k ) $). Conversely, given an 
equivariant (= 
invariant) $k$-form $\theta$ on $P$, and given an infinitesimal 
$k$-simplex $a_0 ,\ldots ,a_k$ in $M$, define 
$$\Theta (a_0 ,\ldots ,a_k ) := \theta (x_0 ,\ldots ,x_k )$$ where
$x_0 ,\ldots ,x_k$ is any infinitesimal $k$-simplex above $a_0 , 
\ldots ,a_k $. The proof that this value does not depend on the choice 
of the $x_i$'s proceeds much like the proof of the well-definedness of 
a connection given a connection-form, in Proposition \ref{one} above: First 
we prove, for fixed $x_0$ above $a_0$, that the value is independent 
of the choice of the remaining $x_i$'s, and this is clear from the 
verticality assumption on $\theta$. Next we prove that changing $x_0$ 
to $x_0 \cdot g$ (and picking $x_1 \cdot g , \ldots ,x_k \cdot g$ for 
the remaining vertices in the new $k$-simplex) does not change the 
value either, and this is clear from equivariance (= invariance).

\medskip

Recall that a $k$-form with values in a group bundle $E\to M$ 
associates to an infinitesimal $k$-simplex $a_0 , ... , a_1$ in 
$M$ an element in the fibre of $E_{a_{0}}$. We are interested 
in the case where $E$ is the gauge group bundle of a groupoid; 
such forms we call {\em gauge forms}, for brevity.    

\begin{prop}
There is  a natural bijective correspondence between 
horizontal equivariant $k$-forms on $P$ with values in $G=P^{-1}P$, and $k$-forms on $M$ 
with values in the gauge group bundle $\mbox{{\bf gauge}}(PP^{-1})$.
\label{3}\end{prop}

{\bf Proof/Construction.} Given a horizontal equivariant 
$k$-form $\theta $ on 	$P$ as above, we construct a gauge 
valued $k$-form $\check{\theta} $ on $M$ by the formula
\begin{equation}
\check{\theta} (a_0 ,\ldots ,a_k ):=(u_0\cdot \theta (u_0 ,\ldots 
,u_k))u_0 ^{-1},
\label{theta-to-alpha}\end{equation}
where
$(u_0 ,\ldots 
,u_k)$ is an arbitrary infinitesimal $k$-simplex mapping to 
the infinitesimal $k$-simplex $(a_0,\ldots ,a_k )$ by $\pi $ 
(such exist, since $\pi $ is a surjective submersion). Note that 
the enumerator and the denominator in the fraction defining 
the value of $\check{\theta} $ are both in the fibre over $x_0$, so 
that the value is an endo-map at $a_0$ in the groupoid 
$PP^{-1}$, thus does belong to the gauge group bundle. --- We 
need to argue that this value does not depend on the choice of 
the infinitesimal simplex $(u_0 ,\ldots u_k )$. We first argue 
that, once $u_0$ is chosen, the choice of the remaining $u_i 
$'s in their respective fibres does not change the value. This 
follows from  (\ref{bb}). To see that the value does not depend on the 
choice of $u_0$: choosing another one amounts to choosing some 
$u_0 \cdot g$, for some  $g$. 
But then we just change $u_1 ,\ldots ,u_k $ by the same $g$; 
this will give the arrow in $PP^{-1}$
$$(u_0\cdot g\cdot \theta(u_0\cdot g,\ldots ,u_k \cdot g))(u_0 
\cdot  g)^{-1}.$$
Now we  calculate using the ``metatheorem'', so we drop 
partentheses and multiplication dots; using the assumed equivariance 
(\ref{ccc}), this expression then yields
$$u_0 g g^{-1} \theta (u_0 ,\ldots ,u_k )gg^{-1} u_0^{-1},$$
which clearly equals the expression in (\ref{theta-to-alpha}).

Conversely, given a gauge valued $k$-form $\alpha $ on $M$, we 
construct a $P^{-1}P$-valued $k$-form $\hat{\alpha}$ on $P$ by putting
 
\begin{equation}
\hat{\alpha} (u_0  ,u_1 ,\ldots ,u_k ):=u_0 ^{-1} (\alpha (a_0 
, a_1 ,\ldots 
,a_k )\cdot u_0 )
\label{alpha-to-theta}\end{equation}
where $a_i $ denotes $\pi (u_i )$. Since, for $i\geq 1$, this expression 
depends on $u_i $ only through $\pi (u_i )=a_i $,  it is 
clear that  (\ref{bb}) holds, so 
the form $\hat{\alpha} $ is horizontal. Also,
$$\hat{\alpha} (u_0 \cdot g,\ldots ,u_k \cdot g)=(u_0\cdot g )^{-1}(\alpha 
(a_0 ,\ldots ,a_k )\cdot (u_0 \cdot g ));$$ by the 
metatheorem, this  
immediately calculates to the expression in 
(\ref{alpha-to-theta}).

Finally, a calculation with the metatheorem again (cancelling 
$u_0 ^{-1}$ with $u_0$) immediately gives that the two processes 
$\theta \mapsto \check{\theta}	$ and
$\alpha \mapsto \hat{\alpha}$ are inverse to each other.

We may summarize the bijection $\alpha \mapsto \hat{\alpha}$ 
from $\mbox{{\bf gauge}}(PP^{-1})$-valued forms on $M$ to 
horizontal equivariant $P^{-1}P$-valued forms on $P$  by the 
formula
\begin{equation}
u_0 \cdot \hat{\alpha}(u_0 ,,,, ,u_k )=(\pi^* \alpha)(u_0 , 
... , u_k )\cdot u_0 . 
\label{hat}\end{equation}
In the case that the group $G=P^{-1}P$ is commutative, we may 
cancel the ``external'' $u_0$'s, and get
$$\hat{\alpha} (u_0 , ... ,u_k ) =(\pi ^* \alpha )(u_0 , ... 
,u_k ),$$
for all infinitesimal $k$-simplices $u_0 , ... ,u_k $. So under the 
identification of gauge forms with $G$-valued forms implied by 
(\ref{comm}), we have  that
\begin{equation}\hat{\alpha} = \pi ^* \alpha .
\label{comm2}\end{equation}

\medskip

Recall that if $\nabla $ and $\nabla _1$ are two connections in 
a groupoid $\Phi \doublearrow M$, we may form a 1-form $\nabla _1 \nabla 
  ^{-1}$ with values in the gauge group bundle; it is 
given by
$$\nabla _1 \nabla  ^{-1}(a,b) = \nabla _1(a,b)\cdot \nabla (b,a).$$

For the case where the groupoid is $PP^{-1}$, we have the 
following Proposition, which we shall not use in the sequel, 
but include for possible future reference:

\begin{prop}Let $P\to M$ be a principal bundle, and let 
$\nabla $ and $\nabla _1$ be two connections in the groupoid 
$PP^{-1}$. Then
$$(\nabla _1 \nabla  ^{-1})\hat{\mbox{}} = \omega _1 \cdot \omega ^{-1}$$
where $\omega $ and $\omega _1$ are the connection forms of $\nabla $ 
and $\nabla _1$, respectively.
\end{prop}
{\bf Proof.} Let $x\sim y$, over $a$ and $b\in M$, 
respectively. Then 
\begin{eqnarray*}(\nabla _1 \nabla  ^{-1})\hat{\mbox{}}(x,y)&=
&x^{-1}(\nabla _1 (a,b)\nabla (b,a) x)\\
&=&x^{-1}\nabla _1 (a,b) y \omega (y,x)\\
&=&x^{-1}x \omega _1 (x,y)\omega (y,x)\\
&=&\omega _1 (x,y)\omega (y,x)\\
& =& (\omega _1 \omega ^{-1})(x,y),
\end{eqnarray*}
using the defining relation (\ref{hat}) for $(-)\hat{\mbox{}}$, and the relation 
(\ref{basic}) for $\nabla$ and $\nabla _1$, respectively.

\section{Curvature versus coboundary} Recall that the {\em 
curvature} of a connection in a groupoid 
$\Phi \doublearrow M$ is the $\mbox{{\bf Gauge}}(\Phi )$-valued 
2-form $R=R_{\nabla }$ given by $$R  (a_0 ,a_1 
,a_2 ) = \nabla (a_0 ,a_1 )\cdot  \nabla (a_1 ,a_2 ) \cdot 
\nabla (a_2 ,a_0 ),$$ and recall that if $\omega$ 
is a 1-form with values in a group $G$, them $d\omega $ is the 
$G$-valued 2-form given by $$d\omega (x_0 ,x_1 , x_2 )= \omega 
(x_0 ,x_1 )\cdot \omega (x_1 ,x_2 ) \cdot \omega (x_2 ,x_0 ).$$

We apply this to the case where $\Phi = PP^{-1}$ and 
$G=P^{-1}P$, for a principal fibre bundle $\pi :P\to M$. Then the 
curvature $R$, which is a $\mbox{{\bf gauge}}(PP^{-1})$
-valued 2-form on $M$, 
gives, by Proposition \ref{3}, rise to a (horizontal and equivariant) $P^{-1}P$-valued 
2-form $\hat{R}$ on $P$.

We then have the following:

\begin{thm}Let $\pi :P\to M$ be a principal fibre bundle, and let 
$\nabla$ be a connection in the groupoid $PP^{-1}$
with connection form $\omega$.  Then we have 
an equality of $P^{-1}P$-valued 2-forms on $P$:
$$\hat{R} = d\omega .$$
In particular, $d\omega$ is horizontal and equivariant. 
\label{GaussBonnet}\end{thm}

The form $\hat{R}=d\omega$ is the {\em curvature form} of the 
connection. See the remark below for comparison with the classical 
formulation.

{\bf Proof.} Let $x,y,z$ be an infinitesimal 2-simplex in $P$, and let 
$a=\pi (x)$, $b= \pi (y)$, and $c= \pi (z)$. We calculate the 
effect of the (left) action of the arrow $R(a,b,c)$ on 
$x$ (note that $R(a,b,c)$ is an endo-arrow at $a$ in 
the groupoid):
\begin{eqnarray*}R (a,b,c)\cdot x &=& \nabla (a,b)\cdot \nabla (b,c) 
\cdot \nabla (c,a) \cdot x\\
 & = &\nabla (a,b)\cdot \nabla (b,c) \cdot z \cdot \omega (z,x)\\
&=&\nabla (a,b) \cdot y \cdot \omega (y,z)\cdot \omega (z,x)\\
&=&x\cdot \omega (x,y)\cdot \omega (y,z) \cdot \omega (z,x)\\
&=&x\cdot d\omega (x,y,z),\end{eqnarray*}
using the defining equations for $R$ and for $dw$ for 
the two outer equality signs, and using (\ref{basic}) three times 
for the middle three ones. This proves the Theorem.

\medskip

{\bf Remark.} By \cite{SDG} I.18, or in more detail, 
\cite{DFVG}), there is a bijective 
correspondence between $G$-valued $k$-forms $\theta$ on a manifold $P$ 
(where $G$ is a Lie group, say $P^{-1}P$), and differential 
$k$-forms $\overline{\theta}$, in the classical sense, with values in the Lie 
algebra $\underline{g}$ of $G$ (i.e.\ multilinear alternating 
maps $TP \times _P ... \times _P TP \to \underline{g}$. Under 
this correspondence, the horizontal equivariant 2-form 
$d\omega$ considered in the Theorem corresponds to the 
classically considered "curvature 2-form" $\Omega $ on $P$, as 
in \cite{Nomizu} II.4, \cite{BC} 5.3, or \cite{CDD} V bis 4, (perhaps modulo 
a factor $\pm 2$, depending on the conventions chosen). This 
is not completely obvious, since $\Omega$ differs from the 
exterior derivative $d\overline{\omega}$ of the classical connection form 
$\overline{\omega}$ by a "correction term" $1/2 [\overline{\omega} 
,\overline{\omega} ]$ involving the Lie Bracket of $\underline{g}$; or, alternatively, 
the curvature form comes about by  modifying $d\overline{\omega}$ by a ``horizontaliztion 
operator'' (this ``modification'' also occurs in the treatment in 
\cite{SIL}).  The fact that 
this ``correction term'' (or the ``modification'') does not come up in  our context 
can be explained by Theorem 5.4 in \cite{DFVG} (or see \cite{SDG} 
Theorem 18.5); here it is proved that the formula $d\omega 
(x,y,z) = \omega (x,y) \cdot \omega (y,z) \cdot \omega (z,x)$ 
already contains this correction term, when translated into 
"classical" Lie algebra valued forms.

\medskip  

For the  case where the group $P^{-1}P$ is commutative, we may 
use the isomorphism (\ref{comm}) to identify $\mbox{\bf{ 
gauge}}(PP^{-1})$-valued forms on $M$ with $P^{-1}P$-valued 
forms on $M$. Also, by Proposition \ref{two}, and the 
horizontality and equivariance of $d\omega$, there is a unique 
$P^{-1}P$-valued 2-form $\Omega = \check{d\omega}$ on $M$ with $\pi ^{*} (\Omega 
)=d\omega $.
We therefore have the following Corollary (notation as above), 
which is essentially what \cite{SIL} call the infinitesimal 
version of Gauss-Bonnet Theorem (for the case where $G=SO(2)$):

\begin{cor}
Assume $P^{-1}P$ is commutative, and let the connection $\nabla $ in 
$PP^{-1}$ have connection form $\omega$. Then the unique $G$-valued 2-form 
$\Omega$ on $M$ with $\pi^{*} \Omega = d\omega$ is $R _{\nabla}$. 
$$R_{\nabla }= 
\Omega .$$
\label{cor}\end{cor}

Let us remark that \cite{SIL} also gives a version of the 
Corollary for the non-commutative case, their Proposition 
6.4.1; this, however, seems not correct. In this sense, our Theorem \ref{GaussBonnet} is 
partly meant as a correction to Prop. 6.4.1, partly a ``translation'' of it 
into the pure multiplicative fibre bundle calculus, which is 
our main concern.


\end{document}